\documentclass{article}
\usepackage{amsfonts}
\newtheorem{theorem}{Theorem}
\newtheorem{remark}[theorem]{Remark}
\newtheorem{corollary}[theorem]{Corollary}
\newtheorem{lemma}[theorem]{Lemma}
\newtheorem{proposition}[theorem]{Proposition}
\newcommand{\ch}{\mbox{\textrm{char\,}\nobreak}}
\newcommand{\epr}{\hfill$\diamondsuit$\smallskip}
\newcommand{\pr}{\textit{Proof}. }
\newcommand{\spanned}{\mbox{\textrm{sp\,}\nobreak}}
\title{2-graded polynomial identities for the Jordan algebra
of the symmetric matrices of order two}
\author{Plamen Koshlukov\thanks{Partially supported by grants from CNPq (Nr. 302651/2008-0), and
 from FAPESP (Nr. 2005/60337-2 and 2010/50347-9)}\\
Department of Mathematics, IMECC, UNICAMP\\
Sergio Buarque de Holanda 651, Cidade Universitaria \\
13083-859 Campinas, SP, Brazil\\
e-mail: \texttt{plamen@ime.unicamp.br},\\
Diogo Diniz P. S. Silva\thanks{Partially supported by PhD grant from FAPESP (Nr. 2007/00447-4)}\\ 
Department of Mathematics, UAME/CCT, UFCG \\
P.O.Box 10044, 58109-970 Campina Grande, PB, Brazil\\
e-mail: \texttt{diogo@dme.ufcg.edu.br}
}

\date{}
\begin{document}
\maketitle
\begin{abstract}
The Jordan algebra of the symmetric matrices of order two over a field
$K$ has two natural gradings by $\mathbb{Z}_2$,
the cyclic group of order 2. We describe the graded polynomial
identities for these two gradings when the base field is infinite and of
characteristic different from 2. We exhibit  bases for these
identities in each of the two cases. In one of the cases we perform
a series of computations in order to reduce the problem to dealing with
associators while in the other case one
employs methods and results from Invariant theory. 
Moreover we extend the latter grading to a
$\mathbb{Z}_2$-grading on $B_n$, the Jordan
algebra of a symmetric bilinear form in a vector space of
dimension $n$ ($n=1$, 2, \dots, $\infty$). We call this
grading the \textsl{scalar} one since its even part consists
only of the scalars. 
As a by-product we obtain finite bases of the $\mathbb{Z}_2$-graded identities for
$B_n$. In fact the
last result describes the weak Jordan polynomial identities for the
pair $(B_n, V_n)$.
\end{abstract}



\section*{Introduction}
Polynomial identities in simple associative algebras over a field have
always been of significant interest to ring theorists. Most of the
research in this area has been done under the assumption of the base
field being infinite, and even stronger, of characteristic 0. In spite
of the extensive research in this area little is known about the
concrete form of the identities satisfied by such algebras. In fact
(assuming the base field $K$ infinite) the identities of the matrix
algebras $M_n(K)$ are known only for $n=2$ and $\ch K\ne 2$
(see for example \cite{razmbook}, \cite{drm2}, \cite{kjam2}). One may
study other types of identities. Thus the trace identities of the
matrix algebra $M_n(K)$, $\ch K=0$, were described independently by Procesi
\cite{procesi}, and by Razmyslov (see for example
\cite{razmbook}). Here we recall that of great importance for Ring
theory have also been the methods developed in the course of studying
the trace identities of $M_n(K)$. Also weak identities were introduced and used
successfully by Razmyslov in describing the identities of associative
and Lie algebras, see for an account and various applications
\cite{razmbook}.

Later on in the 80-ies, a powerful theory developed by Kemer provided
a description of the ideals of identities (also called T-ideals) in
the free associative algebra over a field of characteristic 0. One
finds details concerning Kemer's theory in \cite{kemerbook}. One of
the principal tools in that theory was the
usage of $\mathbb{Z}_2$-graded algebras and their graded
identities. Clearly this provided a strong impulse to the study of graded
identities in associative algebras. The interested reader can look at
\cite{bz} and at \cite{pkijac} and their references for some of the
important results concerning gradings and graded identities.

The latter paper dealt with graded identities in the Lie algebra
$sl_2(K)$, and it is one of the few about the topic.
It is somewhat surprising that graded identities in Jordan algebras
have not been studied in detail yet except for the paper \cite{vasjalg}
where the author described the $\mathbb{Z}_2$-graded identities of the Jordan
superalgebra of a bilinear form.

In this paper we study the $\mathbb{Z}_2$-graded identities for the
Jordan algebra $B_2$ of the symmetric matrices of order two over an
infinite field $K$ of characteristic different from 2. Up to a graded
isomorphism there are two nontrivial $\mathbb{Z}_2$-gradings on
$B_2$. We exhibit finite bases (that is generators) of the
corresponding ideals of graded identities. In one of the cases our
methods are more general than needed and we are able to describe a
basis of the weak Jordan identities of the pairs $(B_n, V_n)$ and
$(B,V)$ where $B_n$ and $B$ stand for the Jordan algebra of a
nondegenerate symmetric bilinear form on the vector spaces $V_n$ and
$V$, respectively, $\dim V_n=n$ and $\dim V=\infty$. It is well known
that the Jordan algebra of the symmetric $2\times 2$
matrices is a Jordan algebra of a nondegenerate symmetric
bilinear form on a vector space of dimension 2 therefore there is nothing wrong with
our notation. Recall that the
Jordan algebras $B_n$ and $B$ are simple and special and their associative
enveloping algebras are the Clifford algebras $C_n$ and $C$,
respectively. The weak (associative) identities of the pairs
$(C_n,V_n)$ and $(C,V)$ were described in \cite{pkfbi} over an
infinite field of characteristic $\ne 2$, see also
\cite{vdpk} for the case $\ch K=0$.

\section{Preliminaries}
Throughout $K$ stands for an infinite field of characteristic
different from 2; all vector spaces and algebras (not necessarily
associative) are considered over
$K$. If $A$ is an associative algebra then $A^+$ stands for the vector
space of $A$ equipped with the Jordan product $a\circ b=(ab+ba)/2$. It
is immediate to see that $A^+$ is a Jordan algebra; the Jordan
algebras of this type and their subalgebras are called
\textsl{special}, otherwise they are \textsl{exceptional}. Let $V$ be
a vector space with a nondegenerate symmetric bilinear form $\langle
u,v\rangle$, and let $B=K\oplus V$. One defines a multiplication
$\circ$ on $B$ as follows. If $\alpha$, $\beta\in K$ and $u$, $v\in V$ then
$(\alpha+u)\circ (\beta+v) = (\alpha\beta + \langle u,v\rangle) +
(\alpha v+\beta u)$. It is not difficult to check that $B$ is a Jordan
algebra. We denote it by $B$ whenever $\dim V=\infty$, and by $B_n$
when $\dim V=n$. 
In fact the above is an abuse of notation
since the algebras $B_n$ and $B$ depend on the form $\langle
u, v\rangle$. Clearly equivalent symmetric bilinear forms
define isomorphic Jordan algebras, and vice versa. Let
us observe that if the field $K$ is algebraically closed
then up to an isomorphism there is only one algebra
$B_n$. Over an arbitrary field one has to interpret $B_n$, respectively $B$,
as a class of Jordan algebras which are not necessarily
isomorphic but are classified by the inequivalent
nondegenerate symmetric bilinear forms on the corresponding
vector space.
In order to keep the notation consistent we shall
denote the vector space in the latter case by $V_n$, that is $\dim
V_n=n$. If $A$ is any algebra one defines the
\textsl{associator} of $a$, $b$, $c\in A$ as $(a,b,c) = (ab)c - a(bc)$
where $ab$ is the product in $A$. Let $J_2$ be the vector space of the
symmetric $2\times 2$ matrices, it is a subalgebra of the Jordan
algebra $M_2(K)^+$. Since for every two traceless matrices $a$, $b\in
M_2$ one has that $a\circ b$ is a scalar multiple of the unit
matrix (and the trace is nondegenerate) one gets 
$J_2\in B_2$. In view of the previous remark we consider
$B_2$ as the class of the Jordan algebras of a nondegenerate
symmetric bilinear form on a vector space $V_2$, $\dim
V_2=2$.

We fix the basis $I=e_{11}+e_{22}$, $a=e_{11}-e_{22}$,
$b=e_{12}+e_{21}$ of $J_2$. Here $e_{ij}$ are the usual matrix
units. One gets immediately $a^2=b^2=I$, $a\circ b=0$.

If $G$ is a group and $A$ is an algebra then $A$ is $G$-graded if
$A=\oplus_{g\in G} A_g$, a direct sum of vector subspaces such that
$A_gA_h\subseteq A_{gh}$ for all $g$, $h\in G$. The elements of $A_g$
are homogeneous of degree $g$. If $a\in A_g$ then we shall denote its
homogeneous degree by $|a|$. We shall consider only
gradings by the \textsl{additive} cyclic group $\mathbb{Z}_2$. So a
graded algebra in this paper means $A=A_0\oplus A_1$ where
$A_iA_j\subseteq A_{i+j}$, $i$, $j=0$, 1, and the latter sum taken
modulo 2. Sometimes, when needed, we shall use upper indices to denote
the graded components of $A$. That is when necessary we shall write $A
= A^{(0)}\oplus A^{(1)}$ instead of $A=A_0\oplus A_1$.

Let $X$ be an infinite countable set, $X=\{x_1,x_2,\ldots\}$, and
denote by $K(X)$ and by $J(X)$ the free (unitary) associative and the
free Jordan algebra freely generated by $X$ over $K$. A polynomial
$f=f(x_1,\ldots, x_n)\in K(X)$ is a \textsl{polynomial identity} (a PI
or an identity) for the associative algebra $A$ if $f(a_1,\ldots,a_n)
= 0$ for every $a_i\in A$. The set $T(A)$ of all identities for $A$ is
an ideal in $K(X)$ that is closed under endomorphisms; such ideals are
called T-ideals. It is easy to show that every T-ideal is the ideal of
identities of certain algebra. In the same manner one defines Jordan
identities and T-ideals in $J(X)$. An identity $f$ is a consequence of
the identity $g$ (or $f$ follows from $g$ as an identity) if $f$ lies
in the T-ideal generated by $g$. Similarly $f$ and $g$ are equivalent
as identities if each of them follows from the other. A set of
identities is called a \textsl{basis} of the T-ideal $I$ if this set
generates $I$ as a T-ideal.

As a general rule finding a basis of the identities of a given algebra
may be an extremely difficult (if not hopeless) problem. Here we point
out that in the associative case, a celebrated result of Kemer states
that every nontrivial T-ideal has a finite basis in characteristic 0
(see \cite{kemerbook}). This had been for more than 30 years the
famous Specht problem. We also remark that except for the matrix
algebras $M_n(K)$, $n\le 2$, and $\ch K\ne 2$, for the Grassmann algebra $E$ and its
tensor square $E\otimes E$ (the latter only when $\ch K=0$) no other
''interesting'' T-ideal has an explicit finite basis known.

The polynomial identities of the
Jordan algebras $B_n$, and $B$ were described by Vasilovsky in
\cite{vasja} under some minor restrictions on the characteristic of
the base field. Namely Vasilovsky found bases of the corresponding
T-ideals. (Recall that Iltyakov in \cite{ilt} dealt with the case of
$B_n$ in characteristic 0.) The structure of the relatively free algebras
$J(X)/T(B_n)$ and $J(X)/T(B)$ was given in \cite {vdjordalg}, \cite{pkjaca}.

Let $A$ be an algebra and $V$ a subspace such that $V$ generates $A$
as an algebra. A polynomial $f(x_1,\ldots,x_n)$ is a \textsl{weak
 identity} for the pair $(A,V)$ if $f(v_1,\ldots, v_n)=0$ for all
$v_i\in V$. Depending on the choice of $A$ and $V$ one defines rules
for consequences of a weak identity. Since we shall deal with a
particular situation we do not treat the most general case here but
instead we refer the reader to \cite{pkfbi} and \cite{vdpkjims} for
further information about weak identities in a general setting.

Let $A=J$ be a Jordan algebra and $V$ a subspace that generates $J$
as an algebra. In this case one speaks of weak Jordan identities. Note
that here the weak identities are polynomials in the free Jordan algebra. We
denote by $T(J,V)$ the set of the weak identities for the pair $(J,V)$
and call it the \textsl{weak T-ideal} of that pair. If $f(x_1,\ldots,
x_n)\in J(X)$ then the weak T-ideal $\langle f\rangle^w$ defined by
$f$ is the ideal of $J(X)$ generated by all polynomials
$f(g_1,\ldots,g_n)$, $g_i\in J(X)$. This determines the rule for taking
consequences of a given polynomial, or set of polynomials, as weak
identities.

We shall need some facts about gradings and graded identities. Let
$X=Y\cup Z$ be a disjoint union of the infinite sets
$Y=\{y_1,y_2,\ldots\}$ and $Z=\{z_1,z_2,\ldots\}$. We define a
$\mathbb{Z}_2$-grading on the free Jordan algebra $J(X)$ as
follows. If $m$ is a monomial then it is of $\mathbb{Z}_2$-degree 0
(that is an even element) if its total degree in the variables $Z$ is
even; otherwise it is of $\mathbb{Z}_2$-degree 1 that is, an odd
element. Put $J(X)_i$ the span of all monomials of
$\mathbb{Z}_2$-degree $i$, $i=0$, 1, then $J(X) = J(X)_0\oplus J(X)_1$
is a grading on $J(X)$. Let $J$ be a graded Jordan algebra,
$J=J_0\oplus J_1$, a polynomial $f(y_1,\ldots,y_m,z_1,\ldots, z_n)$ is
a graded identity for $J$ if $f$ vanishes whenever one substitutes the
variables $y_i$ by any elements of $J_0$ and the $z_i$ by any elements
of $J_1$. Clearly the set $T_2(J)$ of all graded identities of $J$ is
an ideal that is closed under endomorphisms of $J$ that respect its
grading.

We consider unitary algebras only. It is easy to show that the unit
element 1 lies in the even component of the graded algebra. Suppose
$J\in B_n$ is graded, $J=J_0\oplus J_1$. Thus $K\subseteq J_0$. Let the
vector space $J_0$ have a basis consisting of 1 and $v_1$, \dots,
$v_k\in V_n$. If $\alpha+v\in J_1$, $\alpha\in K$, $v\in V_n$ then
$v_i\circ (\alpha+v) \in J_1$ but on the other hand $v_i\circ
(\alpha+v) = \alpha v_i + \langle v_i,v\rangle\in J_0$. Hence
$\alpha=0$ and $\langle v_i,v\rangle =0$. Therefore $J_i$ is a
subspace of $V_n$ and moreover $J_1$ is orthogonal to 
$\spanned(v_1,\ldots,v_k)$, the span of $v_1$, \dots,
$v_k$. The same argument as above transfers
verbatim to $B$. In this way we prove the following proposition.

\begin{proposition}
\label{gradings}
Every $\mathbb{Z}_2$-grading on the Jordan algebras $B_n$ and $B$ is defined by a
splitting of the vector spaces $V_n$, respectively $V$, as a direct sum
of two orthogonal subspaces.
\end{proposition}

In the paper \cite{bsz} the authors described \textsl{all} gradings
on the Jordan algebras $B_n$. So our
Proposition~\ref{gradings} is a particular case of the results in
\cite{bsz}. We give a proof of it here since in our particular case
the proof is rather simple and straightforward without relying on the
technique developed in \cite{bsz}.

\begin{corollary}
\label{gradingsj2}
Up to a graded isomorphism there are two nontrivial $\mathbb{Z}_2$-gradings
on the Jordan algebra $J=J_2$. These are given by $J=J_0\oplus J_1$
where either $J_0=\spanned(I, a)$, $J_1=\spanned(b)$, or
$J_0=\spanned(I)$, $J_1=\spanned(a,b)$.
\end{corollary}

In the next two sections we shall handle the graded identities in each
of the two possibilities for the grading on $J_2$. The reader will
observe that the latter grading has as an even part the 
scalars only. That is why we call the first of the gradings described
in the Corollary the \textsl{nonscalar} and the second the
\textsl{scalar} grading, respectively.

\begin{remark}
\label{z3grading}
The algebra $J_2$ admits nontrivial gradings by groups other
than $\mathbb{Z}_2$. The existence of such gradings may
depend on the ground field $K$. Let us consider for example
$G=\mathbb{Z}_3$, the cyclic group of order 3. Suppose $J_2
= (J_2)_0\oplus (J_2)_1\oplus (J_2)_2$ is a nontrivial $G$-grading on
$J_2$. One concludes, according to \cite[Theorem 1]{bsz}
that $(J_2)_0=K$ and, as $G$ has no elements of order two,
$\dim (J_2)_1=\dim (J_2)_2 = 1$. Furthermore $(J_2)_1=K w_1$
and $(J_2)_2= K w_2$ where $w_1$ and $w_2$ are traceless
symmetric matrices such that $w_1^2=w_2^2=0$ and $w_1\circ
w_2=1$. In other words the span of $w_1$ and $w_2$ is a
hyperbolic plane with respect to the bilinear form. Matrices
with these properties exist if and only if $i=\sqrt{-1}\in
K$. If this is the case one may choose for example 
\[
w_1 =  \left(\begin{array}{cc} i&1\\
    1&-i \end{array}\right), \qquad
w_2 = \frac12 \left(\begin{array}{cc} -i&1\\
    1&i \end{array}\right).
\]
We observe that one may reach the above conclusion by
computing directly, without relying on the general result of
\cite{bsz}. It is easy to show that $\dim (J_2)_0 = 1$ and  
that $(J_2)_1$ and $(J_2)_1$ consist of
traceless matrices. Hence $(J_2)_r\circ (J_2)_r \subseteq
(J_2)_0\cap (J_2)_s = 0$ where $\{r,s\} = \{1,2\}$. Then   
$w_1$ and $w_2$ are traceless singular square zero 
matrices. The remaining computations are easy and immediate. 
\end{remark}

\section{The nonscalar grading}
In this section we fix the grading $B_2=J=J_0\oplus J_1$ on the Jordan
algebra of the symmetric $2\times 2$ matrices given by
$J_0=\spanned(I, a)$, $J_1=\spanned(b)$. First we list some elementary
facts about the graded identities for this grading. Let $T=T_2(J)$ be
the ideal of the graded identities for $J$. In this section we shall
not use the symbol $\circ$ for the product in a Jordan algebra, we
shall use instead the common symbol $\cdot$ for the product and, as a
rule, we shall simply omit it. Recall that $|u|$ stands for the
$\mathbb{Z}_2$-degree of an element $a$ while $\deg a$ is the usual
degree of a homogeneous element in the free algebra.

Denote by $I$ the
ideal of graded identities generated by the polynomials
\begin{eqnarray}
x_1 (x_2 x_3) &-& x_2 (x_1 x_3) \mbox{\textrm{ \ if
  $|x_1|=|x_2|$}}
\label{one}\\
(y_1 y_2,z_1,z_2) &-& (y_1 (y_2,z_1,z_2) + y_2
(y_1,z_1,z_2) - 2z_1 (z_2,y_1,y_2))\label{two}\\
(y_1 y_2,y_3,z_1) &-& (y_1 (y_2,y_3,z_1) + y_2 (y_1,y_3,z_1))\label{three}\\
(z_1 z_2, x_1,x_2) && \label{four}\\
(y_1,y_2,z_1,x,y_3)&-&(y_1,y_3,z_1,x,y_2)\label{five} 
\end{eqnarray}
Here and in what follows the long associators without inner
parentheses will be left normed. That is $(x_1,x_2,x_3,x_4,x_5)$
stands for $((x_1,x_2,x_3),x_4,x_5)$, and so on. Recall that if one
considers the vector space of a Jordan algebra $J$ equipped with the
trilinear composition $(a,b,c)$, $a$, $b$, $c\in J$, then $J$ becomes
a Lie triple system. Passing from Jordan algebras to Lie triple
systems one may prove, see for example \cite[pp.~343, 344]{jacobson}, that
every (long) associator is a linear combination of left normed
ones. We call these left normed associators \textsl{proper} ones.
Also the letter $y$, with or without an index, stands for an even
variable; $z$ with or without index for an odd variable, and $x$ for
any variable (even or odd).

\begin{lemma}
\label{IsubsetT}
The graded identities from (\ref{one}) to (\ref{five}) hold for the Jordan
algebra $J$. In other words $I\subseteq T$.
\end{lemma}

\pr
The proof consists of a straightforward (and easy) computation, so we
omit it.
\epr

Since we consider algebras over infinite fields then every graded
identity is equivalent to a finite collection of multihomogeneous ones
(namely its multihomogeneous components). Therefore we shall consider
only multihomogeneous identities. Our aim in this section is to prove
that in fact $I=T$.

Denote $L=J(X)/I$ where $J(X)$ is the free graded Jordan algebra. We 
shall work in $L$; we shall keep the same notation for the images 
of the variables $y$ and $z$ in $L$. The graded identity (\ref{four})
implies that the elements $z_iz_j$ lie in the associative centre of
$L$. Therefore the subalgebra of $L$ generated by all $z_iz_j$ 
is associative.

Since the ideal $I$ is homogeneous in the grading it follows that the
algebra $L$ is graded, the grading on it being induced by that on
$J(X)$, so $L=L_0\oplus L_1$.

\begin{proposition}
\label{l0l1}
The subalgebra $L_0$ of $L$ is associative. Also the subalgebra of $L$
generated by $L_1$ is associative.
\end{proposition}

\pr
The graded identity (\ref{one}) implies that $(y_1,y_2,y_3) = 0$ in
$L$. Hence $L_0$ is associative. Now let $w\in [L_1]$, the subalgebra
in $L$ generated by $L_1$. Suppose $w$ is a monomial, $\deg
w=n$. First we prove that $w=(u) \hat{z_i}$ for some $z_i$ and some
monomial $u\in L_C$, the associative centre of $L$. Here the hat over
the variable $z_i$ means that it may be missing (in which case we
consider it as 1). We induct on $n$. If
$n=1$ then $w=z$ and $u=1$. Suppose that $n>1$ and furthermore that
the statement is true for monomials of degree up to $n-1$. Then
$w=(w_1)\ldots (w_k)$ for some  $w_i\in L_1$ (there might be any
distribution of the parentheses). By the induction we may assume $w_i =
(u_i)\widehat{z_{t_i}}$. Since the $u_i$ lie in the associative centre we
reduce the case to $w=(u)(w')$ where $u\in L_C$ and $w'$ is a product
of some variables $z$ with some distribution of parentheses. But by the
identity (\ref{four}) the product of every two variables $z$ lives in the
associative centre $L_C$, so it may be transferred to $u$. Hence if
the degree of $w'$ is even then $w\in L_C$ and if it is odd then
$w=(u')z$ where $u'\in L_C$.

Now let $w_i = (u_i) \widehat{z_{t_i}}$, $i=1$, 2, 3, be three elements of
$[L_1]$. Then we obtain $(w_1,w_2,w_3) = (u_1u_2u_3) (\widehat{z_{t_1}},
\widehat{z_{t_2}}, \widehat{z_{t_3}}) = 0$ in $L$ since by the
identity (\ref{one}) one has $(z_1,z_2,z_3) = 0$. If some of the
$\widehat{z_{t_i}}$ is missing it is substituted by 1, and then the
associator vanishes.
\epr

Let $\Omega\subseteq L$ be the smallest subset of $L$ with the following
property. If $f_1$, $f_2$, $f_3\in \Omega\cup X$ then
$(f_1,f_2,f_3)\in \Omega$. The elements of $\Omega$ are called
\textsl{associators}.

Denote by $J(X)^{(n_1,\ldots,n_k)}$ the multihomogeneous component of
$J(X)$ of all polynomials that are homogeneous of degree $n_i$ in
$x_i$, $1\le i\le k$, and of degree 0 in each of the remaining
variables. One defines similarly $L^{(n_1,\ldots,n_k)}$. We choose the
subset $\Omega_0$ of $\Omega$ as follows. If $\Omega\cap L^{(n_1,\ldots,n_k)}\ne 0$ then we pick one arbitrary nonzero element of this intersection to be in $\Omega_0$, and there are no other elements
in $\Omega_0$.

Now we define a set $A\subseteq L$. It consists of the elements of the
following four types.
\begin{itemize}
\item[(i)]
$(y_{i_1} \ldots y_{i_k}) (z_{j_1}\ldots z_{j_t})$;
\item[(ii)]
$(y_{i_1} \ldots y_{i_k}) u^1$;
\item[(iii)]
$(y_{i_1} \ldots y_{i_k}) (z_{j_1} u^1)$;
\item[(iv)]
$(y_{i_1} \ldots y_{i_k}) u^0$.
\end{itemize}
Here $k$, $t\ge 0$, and $u^i\in\Omega_0$ is an associator, $|u^i| =
i$, $i=0$, 1. We require further that $\deg u^i\ge 3$, that
is the $u^i$ are associators but not
variables. Proposition~\ref{l0l1} allows us to omit the
parentheses in the above expressions.

Let $S=\spanned(A)$ be the span of $A$ in $L$. We shall show that
$L=S$. In order to do that first we prove that the elements of
$\Omega$ are equal in $L$, up to a sign, to those of
$\Omega_0$. We split the proofs of these claims into several lemmas.

\begin{lemma}
\label{morepi}
The following polynomials  lie in $I$.
\begin{itemize}
\item[(a)]
$(x_1,x_2,x_3)$, $|x_1| = |x_3|$;
\item[(b)]
$(y_1z_1, y_2,y_3) - y_1(z_1,y_2,y_3)$;
\item[(c)]
$(z_1,y_1,\ldots, y_{2k}) - (z_1, y_{\sigma(1)}, \ldots,
y_{\sigma(2k)})$ for every permutation $\sigma$ in the symmetric group
$S_{2k}$;
\item[(d)]
$(z_1,y_1,\ldots,y_{2k},z_2,y_{2k+1}) - (z_{\tau(1)},
y_{\sigma(1)},\ldots, y_{\sigma(2k)}, z_{\tau(2)}, y_{\sigma(2k+1)})$
for all $\sigma\in S_{2k+1}$ and $\tau\in S_2$.
\end{itemize}
\end{lemma}

\pr
Since $(x_1,x_2,x_3) = (x_1x_2)x_3 - x_1(x_2x_3) = (x_3x_2)x_1 -
x_1(x_2x_3) = 0$ in $L$, by (\ref{one}) we obtain
(a). Similarly by (\ref{one}) we have
$((y_1z_1)y_2)y_3 = ((y_2z_1)y_3)y_1$ and $(y_1z_1)(y_2y_3) =
((y_2y_3)z_1)y_1$. Therefore
\begin{eqnarray*}
(y_1z_1,y_2,y_3) &=& ((y_1z_1)y_2)y_3 - (y_1z_1)(y_2y_3) \\
&=&((y_2z_1)y_3)y_1 - ((y_2y_3)z_1)y_1 = y_1(z_1,y_2,y_3)
\end{eqnarray*}
which proves (b). It is clear that (c) follows from
the graded identity (\ref{five}).

In order to prove (d) it suffices to
consider first the case $\sigma=1$, the identity permutation, and then
apply the graded identity (\ref{five}). Suppose first $k=1$. It
follows from (\ref{one}) that $((y_1y_2)z_1)z_2 = ((y_1y_2)z_2)z_1$
and that
\[
(y_1(y_2z_1))z_2 = (y_1z_2)(y_2z_1) = (y_2(y_1z_2))z_1 =
(y_1(y_2z_2))z_1.
\]
Therefore we obtain
\begin{equation}
\label{firstadd}
(y_1,y_2,z_1)z_2 = (y_1,y_2,z_2)z_1
\end{equation} 
It follows from the latter identity, together with (b) and (a), that
\begin{equation}
\label{secondadd} 
(y_1,y_2,z_1)(z_2y_3) = (y_1,y_2, z_2y_3)z_1 = (y_3(y_1,y_2,z_2))z_1 =
(y_1,y_2,z_2)(y_3z_1)
\end{equation}
Now write $(z_1,y_1,y_2,z_2,y_3) = ((z_1,y_1,y_2) z_2)y_3 -
(z_1,y_1,y_2)(z_2y_3)$. Then according to the graded
identity (\ref{firstadd}) we transpose $z_1$ and $z_2$ in
the first summand on the right-hand side. Moreover according
to (\ref{secondadd}) we do the same in the second summand,
and this finishes the case $k=1$. 

Let $k>1$, and suppose that for every integer $\le k-1$ the
polynomial of (d) lies in $I$. In order to settle this case it suffices to check 
that the following equalities hold in $L$.
\begin{eqnarray*}
(z_1,y_1,\ldots,y_{2k})z_2& =& (z_2,y_1,\ldots,y_{2k})z_1;\\
(z_1,y_1,\ldots,y_{2k})(y_{2k+1}z_2) &=& (z_2,y_1,\ldots,y_{2k})(y_{2k+1}z_1).
\end{eqnarray*}
In order to prove the former is a graded identity in $L$ we
observe that by (\ref{firstadd})
\[ 
(z_1,y_1,\ldots,y_{2k})z_2 = (z_2,y_{2k-1},
y_{2k})(z_1,y_1,\ldots,y_{2k-2}).
\]
Now we induct once again on $k$, supposing that
\[
(z_1,y_1,\ldots,y_{2k-2})z_2  = (z_2,y_1,\ldots,y_{2k-2})z_1.
\]
Then we obtain
\[ 
(z_2,y_{2k-1}, y_{2k})(z_1,y_1,\ldots,y_{2k-2}) =
z_1(z_2,y_{2k-1},y_{2k}, y_1,\ldots,y_{2k-2}).
\] 
Now we apply (c) and get the required equality. The latter equality is
verified in a similar way. By (\ref{secondadd}), and by (b) and (a)
(and once again by induction) we have 
\begin{eqnarray*}
(z_1,y_1,\ldots,y_{2k}) (y_{2k+1}z_2) &=&
(z_2,y_{2k-1},y_{2k})(y_{2k+1}(z_1,y_1,\ldots,y_{2k-2}))\\
&=& ((z_2,y_{2k-1},y_{2k})y_{2k+1})(z_1,y_1,\ldots,y_{2k-2})\\
&=& (z_1y_{2k+1})(z_2,y_{2k-1},y_{2k},y_1,\ldots,y_{2k-2}).
\end{eqnarray*}
Now apply (d) to get the required equality.

Thus both equalities hold in $L$; the identity (d) follows easily from them. 
\epr

We shall need some more graded identities in $L$.

\begin{lemma}
\label{stillmorepi}
The following polynomials are graded identities for $L$.
\begin{itemize}
\item[(i)]
$(y_1,z_2,(y_2z_1)) - (y_2(y_1,z_1,z_2) + z_1(y_1,y_2,z_2))$;
\item[(ii)]
$z_1(z_2,z_3,y_1)$;
\item[(iii)]
$(z_1z_2)(z_3,x,y_1) - (z_1,z_2,y_1,x,z_3)$; 
\item[(iv)]
$(y_1,z_1,z_2)(y_2,z_3,z_4) - z_1(y_1,z_2,z_3,y_2,z_4)$;
\item[(v)]
$(y_1,y_2,z_1)(y_3,y_4,z_2) - z_1(z_2,y_1,y_2,y_3,y_4)$.
\end{itemize}
In other words the above polynomials lie in $I$.
\end{lemma}

\pr
In order to prove (i) observe that $(y_1,z_2,(y_2z_1)) =
(y_1z_2)(y_2z_1) - y_1(z_2(y_2z_1))$, and that
\[
y_1(z_2(y_2z_1)) = y_1(y_2(z_1z_2) + (y_2,z_1,z_2)) = (y_1y_2)(z_1z_2)
+ y_1(y_2,z_1,z_2).
\]
Also we have the following equality in $L$
\begin{eqnarray*}
(y_1z_2)(y_2z_1) &=& z_1(y_2(y_1z_2)) = z_1((y_1y_2)z_2 -
(y_2,y_1,z_2)) \\
&=& (y_1y_2)(z_1z_2) + (y_1y_2,z_2,z_1) - z_1(y_1,y_2,z_2).
\end{eqnarray*}
Subtracting the last two identities and applying the graded identity
(\ref{two}) we obtain (i).

Similarly (ii) follows from (\ref{one}) and (\ref{four}) since
\[
z_1((z_2z_3)y_1) = (z_2z_3)(z_1y_1) = ((z_2z_3)z_1)y_1 =
((z_1z_2)z_3)y_1
\]
and moreover $z_1(z_2(z_3y_1))= (z_1z_2)(z_3y_1) =
((z_1z_2)z_3)y_1$.

The graded identity (iii) holds since $z_1z_2$ is in the associative
centre of $L$ and furthermore, by (\ref{four}) one has
$(z_1(z_2y_1),x,z_3)=0$ in $L$.

Now we deduce (iv). It follows from the graded identities (\ref{four})
and (\ref{one}) that
\begin{eqnarray*}
((y_1z_1)z_2)(y_2,z_3,z_4) &=& (y_2,z_3((y_1z_1)z_2),z_4) =
(y_2,(z_1z_2)(y_1z_3),z_4) \\
&=& (z_1z_2)(y_2,(y_1z_3),z_4).
\end{eqnarray*}
On the other hand $((y_1z_1)z_2)(y_2,z_3,z_4) =
(z_1z_2)(y_1(y_2,z_3,z_4) + z_3(y_2,y_1,z_4))$ by (i). Since $z_1z_2$
is in the (associative) centre we have
\[
((y_1z_1)z_2)(y_2,z_3,z_4) = ((z_1z_2)y_1) (y_2,z_3,z_4) +
((z_1z_2)z_3)(y_2,y_1,z_4).
\]
Hence $(y_1,z_1,z_2)(y_2,z_3,z_4) =
((z_1z_2)z_3)(y_2,y_1,z_4)$. By (\ref{four}) the elements
$z_1z_2$ and $(y_1z_2)z_3$ are in the centre. Thus
$z_1(y_1,z_2,z_3,y_2,z_4) = (z_1(z_2z_3))(y_1,y_2,z_4)$. Now
$z_1(z_2z_3) = (z_1z_2)z_3$ hence $(y_1,z_1,z_2)(y_2,z_3,z_4) =
z_1(y_1,z_2,z_3,y_2,z_4)$.

It remains to prove (v). In the proof of Lemma~\ref{morepi} (d) we
showed that $((y_1y_2)z_1)z_2 = ((y_1y_2)z_2)z_1$. Therefore, by
(\ref{one})
\begin{eqnarray*}
(y_1,y_2,z_1)(y_3,y_4,z_2) &=& (y_1,y_2,(y_3,y_4,z_2))z_1 =
-(y_3,y_4,z_2,y_1,y_2)z_1 \\
&=& (z_2,y_3,y_4,y_1,y_2)z_1.
\end{eqnarray*}
Now we apply the graded identity from Lemma~\ref{morepi} (c) and order
the variables $y$ in the last associator.
\epr

The following proposition shows that the choice of the elements of
$\Omega_0$ can be really arbitrary.

\begin{proposition}
\label{omega0}
Let $u_1$ and $u_2$ be two nonzero associators in $L$ of the same
multidegree. Then $u_1=\pm u_2$.
\end{proposition}

\pr
Let $u\in \Omega$. We start with the following easy remark. Let us 
substitute every even variable of $u$ by the matrix $a$, and 
every odd variable of $u$ by the matrix $b$. Then the
evaluation of $u$ on the algebra $J$ will be $\pm a$ whenever $|u|=0$,
and $\pm b$ whenever $|u|=1$. This statement is checked by an obvious
induction on $\deg u$. If $\deg u=3$ then we can have  
$(a,a,b) -(b,a,a) = b$ and $(b,b,a) = -(a,b,b) = a$. If
$u=(u_1,u_2,u_3)$ then the $u_i$ are associators of lower degree than
$u$, and we apply the induction. 

As we observed earlier, every associator $u\in \Omega$ is a linear
combination of proper ones. We shall prove the proposition for proper
(that is left normed) associators. First we shall show that every such
associator $u$ can be 
written as $u=(z_{i_1}\ldots z_{i_{2m}}) u_t$ where $t=0$, 1, and $u_0
= (z_{i_{2m+1}}, y_{j_1}, \ldots, y_{j_{2k}})$, while $u_1 =
(z_{i_{2m+1}}, y_{j_1}, \ldots, y_{j_{2k}}, z_{i_{2m+2}},
y_{j_{2k+1}})$. Here $m\ge 0$.

We shall induct on the total degree $n$ of $u$ in the variables $z$,
and furthermore on $\ell$ where the total degree $\deg u=2\ell+1$. If $n=0$ there are
no such nonzero associators. Suppose $n=1$. If $u=(x_1,x_2,x_3,
\ldots)$ then exactly one of $x_1$ and $x_3$ is an odd variable, hence
we may assume (up to a sign) it is $x_1$. Then apply
Lemma~\ref{morepi} (c), and the result holds for every $\ell$.

Suppose $n=2$. If $\ell=1$ then $\deg u=3$, and $u=(z_1,z_2,
y_1)=(z_2,z_1,y_1)$. Take $\ell\ge 2$. One cannot have
$u=(z_1,z_2,y_1,\ldots)$ since the dots would stand for even variables
(at least two), and $u=0$. Thus $u=(z_1, y_1, \ldots, y_p, z_2,
y_{p+1},\ldots)$ with $p\ge 1$. (The indices of the variables may be
permuted but we use this simpler notation.) Also the integer $p$ is even since
otherwise $u=0$. Moreover the rightmost dots stand for even
variables. Since the associator $(z_1,y_1,\ldots,y_p,z_2)$ is even (in
the grading) and we have $n=2$ odd variables the rightmost dots are
actually missing. In this way $u=(z_1, y_1, \ldots, y_p, z_2, y_{p+1})$
where $p$ is even, and this case is done for every $\ell$ by
Lemma~\ref{morepi} (d).

Now let $n\ge 3$. We shall prove in this case that $u = (z_{i_1}\ldots
z_{i_{2m}}) u'$ where $u'$ contains 1 or 2 odd variables. Write
$u=(A_1,x_1,x_2)$ where $A_1$ is an associator. If $A_1$ contains at
least 3 odd variables then by the induction ($\deg A_1 = \deg u-2$) we
have $A_1=(z_{i_1}z_{i_2})A_2$. Here $A_2$ is some proper associator
(or a linear combination of such). On the other hand $z_1z_2$ is in the
associative centre therefore $u=(z_{i_1}z_{i_2}) (A_2,x_1,x_2)$. In
this way we can apply the induction to $(A_2,x_1,x_2)$.

If, on the other hand, $A_1$ has only one odd variable then $|A_1|=1$
and $|x_2|=0$ since $u\ne 0$. Thus we get $n=1$ or 2, but this is
impossible due to $n\ge 3$.  Therefore we have to deal with the case
when $A_1$ has exactly two odd variables. Once again by induction
($n=2$) we can assume $A_1=(z_1,y_1,\ldots,y_{2k-2}, z_2,y_{2k-1})$
up to some permutation of the even and, separately, of the odd
variables. Then $|A_1|=0$ and therefore $x_2$ must be some $z$, and
since $n=3$, $x_1$ is some $y$. In order to simplify the notation we
write $u=(A_1,y,z)$ and $A_1=(A_2, z', y')$ where $A_2$ is an
associator, $|A_2|=1$ and $A_2$ contains exactly one odd
variable. Therefore $u=(A_2,z',y',y,z)$. Now we apply the identity
(iii) from Lemma~\ref{stillmorepi} and then use the fact that the
algebra generated by $L_1$ is associative (and commutative) in order
to get
\begin{eqnarray*}
u &=& \pm (A_2z')(y',y,z) = \pm (A_2z')(z,y',y) \\
&=&\pm ((A_2z')z,y',y) =\pm ((z'z)A_2,y',y) = \pm (z'z)(A_2,y',y).
\end{eqnarray*}
It is clear from the last argument that we can permute the variables
$z$ at will; for the variables $y$ it follows from
Lemma~\ref{morepi}.

Now let $u$ and $w$ be two associators (not necessarily proper
ones), of the same multihomogeneous degree. We write each of them as a
linear combination of proper ones, and then apply the 
results proved above for the proper associators. Therefore $u$ and $w$ differ only
by a scalar multiple. But this multiple must be either 1 or $-1$ due
to the remark made at the beginning of the proof. 
\epr

Another consequence of Lemma~\ref{stillmorepi} and the proof of
Proposition~\ref{omega0} is the following.

\begin{corollary}
\label{substinA}
Let $u\in A$. If we substitute any variable $x$ of $u$ by an
associator $w$ such that $|x|=|w|$ then we get a linear combination of
elements of $A$.
\end{corollary}

\pr
The algebra $L_0$ is associative and commutative, and the same holds
for the subalgebra of $L$ generated by $L_1$. If, in some
substitution, there appears an element of the type $z_iz_j$ it can be
''eaten'' by the associator (or by the element $z_{j_1}\ldots z_{j_t}$
in the case of elements of type (i)) in the definition of $A$. We
finish the proof of the corollary by means of a case-by-case analysis,
substituting in each of the elements of $A$, associators for some
variable. These cases are straightforward; one needs also the graded
identities from Lemma~\ref{stillmorepi}. For example, if we substitute
the variable $y$ for $(y_2,z_3,z_4)$ in $(y_1,z_1,z_2)y$, by
Lemma~\ref{stillmorepi} we conclude that we get
$z_1(y_1,z_2,z_3,y_2,z_4)$.
\epr

Recall that we denote by $S$ the span of the set $A$ (defined just
before Lemma~\ref{morepi}). The following lemmas assure that certain
elements belong to $S$.

\begin{lemma}
\label{yyyyz}
The polynomial $N=((y_1\ldots y_k),y,z)$ lies in the span of the
elements of type {\rm(ii)} of the set $S$.
\end{lemma}

\pr
Let $S''$ be the set of the elements of type (ii) of the same
multidegree as $N$. We shall show that $N$ lies in the span $V$ of
$S''$. We induct on $k$. If $k=1$ we have nothing to prove. Then
write, using the graded identity (\ref{three})
\[
N=(y_1\ldots y_{k-1})(y_k, y, z) + y_k((y_1\ldots y_{k-1}), y, z).
\]
The element $(y_1\ldots y_{k-1})(y_k, y, z)$ is clearly of
the type (ii). 
In order to prove that $y_k((y_1\ldots y_{k-1}), y, z)\in V$ we
apply the inductive assumption to the element $((y_1\ldots
y_{k-1}),y,z)$.  
Hence it suffices to prove that all elements of the type 
\[
y((y_1\ldots y_p)(z, y_{p+1},\ldots, y_q)), \quad p<k, \quad q-p\equiv
0\pmod{2}
\]
are linear combinations of elements of type (ii). But the latter element
equals
\[
(y (y_1\ldots y_p))(z, y_{p+1},\ldots, y_q) - (y ,(y_1\ldots y_p), (z,
y_{p+1},\ldots, y_q)).
\]
The first summand is of type (ii). The second summand equals, up to a
sign, $((y_1\ldots y_p), y, (z, y_{p+1},\ldots, y_q))$. First consider
the element $((y_1\ldots y_p), y, z)$. Applying to it the
graded identity (\ref{three}) several times we obtain a linear combination of
elements of the type (ii) from the set $A$. Now according to the
previous Corollary~\ref{substinA}, if we substitute an associator for
a variable in an element of $A$, we get once again elements of $A$ as
long as the $\mathbb{Z}_2$-degree is preserved. This finishes the
proof.
\epr

\begin{lemma}
\label{yyyzz}
The polynomial $N=((y_1\ldots y_k), z_1,z_2)$ lies in $S$.
\end{lemma}

\pr
As in the previous lemma we induct on $k$. We shall show that $N$ is
in the span $V$ of the elements of types (iii) and (iv). The base of
the induction $k=1$ is obvious. It follows from the graded identity
(\ref{four}) that
\[
N=(y_1\ldots y_{k-1}) (y_k,z_1,z_2) + y_k((y_1\ldots y_{k-1}),
z_1,z_2) - 2z_1 (z_2, (y_1\ldots y_{k-1}), y_k).
\]
The first summand from the right is an element of $S$ (of type
(iv)). By the induction we can assume that $((y_1\ldots y_{k-1}),
z_1,z_2)\in V$ is a linear combination of elements of types (iii) and
(iv). Moreover the elements of types (iii) and (iv) are products of
even elements; therefore they lie in the associative algebra $L_0$. Thus
multiplying these by $y_k$ yields once again elements of the same
types. It remains to prove that the last summand lies in $V$. Applying
Lemma~\ref{yyyyz} to $(z_2, (y_1\ldots y_{k-1}), y_k)$ we write it as
a linear combination of elements of type (ii). Therefore it suffices
to prove that elements of the form $M=z_1((y_1\ldots y_n)(z_2,y_{n+1},
\ldots, y_m))$, $n<k$, $m-n$ even, are in $V$. But we have
\[
M=(y_1\ldots y_n)(z_1(z_2,y_{n+1}, \ldots, y_m)) - (z_1,(z_2,y_{n+1},
\ldots, y_m), (y_1,\ldots, y_n)).
\]
Here the first summand on the right is of type (iv). The second is
also in $V$ due to the inductive assumption combined with
Corollary~\ref{substinA}.
\epr

\begin{lemma}
\label{sz}
If $s\in S$ then $sz\in S$.
\end{lemma}

\pr
First we notice that the elements of $A$ of the types (ii), (iii),
(iv) can be obtained by the elements of type (i) after the
substitution of a variable $x$ by an associator $u$ such that
$|u|=|x|$. Thus the lemma will follow from Corollary~\ref{substinA}
if we prove that $((y_1\ldots y_k)(z_1\ldots z_p))z\in S$. If the
number $p$ is even then $z_1\ldots z_p$ lies in the associative centre
of $L$ and hence it will ''eat'' the variable $z$ and we get an
element of type (i). So suppose $p$ is odd. Then the product
$z_2\ldots z_p$ is central and it suffices to show that the element
$R=((y_1\ldots y_k)z_1)z\in S$. One sees easily that
\[
R= ((y_1\ldots y_k), z_1,z) + (y_1\ldots y_k)(z_1z).
\]
But the first summand on the right lies in $S$ due to
Lemma~\ref{yyyzz}, while the second is already of type (i).
\epr

\begin{lemma}
\label{yyz}
The element $N=((y_1\ldots y_k), (y_{k+1}\ldots y_n),z)\in S$.
\end{lemma}

\pr
We induct on $n-k$. If $n-k=1$ this is Lemma~\ref{yyyyz}. Suppose
$n-k>1$. The lemma will follow from the following claim. If we
substitute in $(y_1\ldots y_r)(z, y_{r+1}, \ldots, y_s)$ a variable $y$
by a product of $n-k$ variables $y$ the resulting expression lies in
$S$. The claim clearly holds (for products of any length) if we
substitute some of the variables $y_1$ to $y_r$. If we substitute some
of $y_{r+1}$, \dots, $y_s$, then we first apply Lemma~\ref{yyyyz}, and
then the induction.
\epr

\begin{lemma}
\label{last}
The element $R=(y_1\ldots y_k)((y_{k+1}\ldots y_n)z)$ lies in $S$.
\end{lemma}

\pr
The proof follows from Lemma~\ref{yyz} by using an argument similar to
that of the proof of Lemma~\ref{sz}.
\epr

\begin{proposition}
\label{AspansL}
The set $A$ spans the relatively free graded algebra $L$.
\end{proposition}

\pr
First we claim that the elements $R= ((y_1\ldots y_p)z_1)
((y_{p+1} \ldots y_q)z_2)\in S$. Indeed since $(z_1,y,z_2)$
is a graded identity we have that $R$ can be written as $R = z_1 ((y_1\ldots
y_p)(y_{p+1}\ldots y_q)z_2)$. Therefore our claim follows
from Lemma~\ref{sz} and from Lemma~\ref{last}. Moreover the
product of an even number of elements of $L_1$ lies in the
associative centre of $L$. Hence the fact that $R\in S$,
together with Corollary~\ref{substinA} imply that the
product of two elements from $A$ lies in $S=\spanned
(A)$. Thus $S$ is a subalgebra of $L$. Since $X\subseteq S$
by definition, and $X$ generates $L$ as an algebra we obtain
$S=L$ as required.
\epr

Let $u_1$, $u_2\in A$. We shall call $u_1$ and $u_2$
\textsl{similar} if
\[
u_1 = (y_{i_1}^{n_1}\ldots y_{i_k}^{n_k}) a_1, \quad
u_2 = (y_{i_1}^{n_1}\ldots y_{i_k}^{n_k}) a_2.
\]
Here the $a_i$, $i=1$, 2, are of the form $(z_{j_1}\ldots
z_{j_p}) w_i$, $p\ge 0$, and $w_i$ are associators. Note that
we do not require that $a_1=a_2$. In other words $u_1$ and
$u_2$ are similar if the even variables that appear in them
outside the associators, are the same (counting the
multihomogeneous degrees).

We have gathered all necessary information for the main
result in this section.

\begin{theorem}
\label{mainthm}
Let $K$ be an infinite field, $\ch K\ne 2$. Then the ideal
$T$ of the graded identities for the Jordan algebra $J$ of
the symmetric $2\times 2$ matrices is generated (as a graded
T-ideal) by the identities (\ref{one}) to (\ref{five}). In
other words $T=I$.
\end{theorem}

\pr
We split the proof into three steps.

\textit{Claim} 1. Let $u_1$, \dots, $u_n$ be elements of $A$
having the same multidegree. Suppose no two of them are
similar and that $\sum \alpha_i u_i\in T$ is a graded identity for
$J$ where $\alpha_i\in K$ are scalars. Then all $\alpha_i=0$.

Let $u_i = c_i a_i$ where $a_i$ are as in the definition of
similarity, and the $c_i$ are products of even
variables. Then $c_i\ne c_j$ whenever $i\ne j$. Since $L_0$
is associative (and commutative) we can assume that the
variables $y$ in each $c_i$ are written in ascending
order.

Let $\sum \alpha_i u_i = f(y_1,\ldots, y_p, z_1,\ldots,
z_q)$. Suppose further $f\ne 0$. Define
\[
g(y_1,\ldots, y_p,z_1,\ldots,z_q) = f(y_1+1,\ldots,y_p,
z_1,\ldots, z_q).
\]
The polynomial $g$ is a graded identity for the Jordan
algebra $J$. We draw the reader's attention that $g$ is not
multihomogeneous. Since the base field is infinite then all
its multihomogeneous components are also graded identities
for $J$. One of its homogeneous components is exactly
$f$. Take the homogeneous component $h$ of $g$ that is
nonzero and of the least degree in $y_1$. (That is we take
for $h$ the nonzero polynomial obtained from $f$ after substituting
the largest possible number of variables $y_1$ by 1.) The
polynomial $h$ is obtained from $f$ by means of the
following procedure. First we take the sum of all $\alpha_i
c_iu_i$ where the degree of $y_1$ in $c_i$ is the largest
possible, and discard the remaining summands. Then we
substitute in these summands, all entries of $y_1$ in $c_i$
by 1 (and keep the entries of $y_1$ in the
associators.). This gives exactly $h$ since whenever 1 is
substituted in an associator, the associator vanishes. Now the polynomial
$h$ does not contain $y_1$ outside associators. By repeating
the above argument to $h(y_1,y_2+1,y_3,\ldots,y_p,z_1,\ldots, z_q)$ we shall get a nonzero polynomial that
does not contain $y_2$ outside associators, and so
on. Finally we shall get a non-zero polynomial $f_1$ that does not contain
any variable $y_i$ outside its associators. Clearly $f_1\in T$ since
$f\in T$. But $f_1$ is obtained by $f$ by removing some of the
summands and discarding the $c_i$ parts of the remaining summands. Now
as the $c_1$, \dots, $c_n$ are pairwise distinct we get that there is
only one $a_i$ in $f_1$. That is $f_1=\alpha_i a_i$ for some $i$.
On the other hand if $\alpha_i a_i\in T$ this means it must be a
graded identity for $J$. But this is possible only if $\alpha_i=0$ in
which case $f_1=0$, and $\alpha_i c_i a_i$ does not participate in
$f$. Then we repeat the above procedure to $f$ (having discarded
$\alpha_i c_i a_i$) and we  continue by induction.

\medskip

\textit{Claim} 2. The set $A$ is linearly independent modulo
the graded ideal $T$.

It follows from Claim 1 that it suffices to consider only
the elements of $A$ where all variables $y$ appear in the
associators only. Thus we have to show that the elements
$z_{j_1}\ldots z_{j_t}$, $u^1$, $z_{j_1} u^1$, $u^0$, are
independent. Note that that the $u^i$ are associators,
not variables. Then these elements are of pairwise different
multidegrees, and cannot be linearly dependent. Our claim
follows.

\medskip

\textit{Claim} 3. The inclusion $T\subseteq I$ holds.

Let $f\in T$ be a multihomogeneous polynomial. Since
$I\subseteq T$, by Proposition~\ref{AspansL} it follows that
$f\equiv \sum \alpha_i u_i\pmod{I}$. Here $\alpha_i\in K$
and $u_i\in A$. But Claims 1 and 2 yield that all
$\alpha_i=0$, and $f\in I$. The claim is proved.

\medskip

In order to finish the proof of the theorem it suffices to
recall that $T\subseteq I$ and $I\subseteq T$.
\epr

\section{The scalar grading}
In the scalar grading on $J$, as we commented before, the component
$J_0$  is the one dimensional span of the unit matrix $I$, and hence
we identify $J_0=K$. It turns out that the scalar grading is somewhat ''easier'' to resolve even in a more general
situation. Here we shall describe the graded identities of the Jordan
algebras $B$ and $B_n$ of a nondegenerate symmetric bilinear form on 
the vector spaces $V$ and $V_n$, respectively. Here $\dim V=\infty$,
$\dim V_n=n$. We have $B^{(0)} = K$, and $B^{(1)}=V$, respectively
$B_n^{(0)} = K$, $B_n^{(1)}=V_n$. (We
shall use upper indices for the grading in order not to
confuse them with the corresponding algebras $B_0$ and $B_1$.) We keep
the notation $X$ for the variables in the free
Jordan algebra, $X=Y\cup Z$ where $Y$ are the even and $Z$ are the odd
variables. We have the following  graded identity in $B$.
\begin{equation}
\label{mainpiscalar}
(y,x_1,x_2) = 0
\end{equation}

Its validity in $B$ is immediate: the even elements are scalars. Note
that it follows from (\ref{mainpiscalar}) that $(z_1z_2, x_1,x_2)=0$. 

We denote in this section by $I$ the ideal of graded identities
defined by the polynomial (\ref{mainpiscalar}), and by $L=J(X)/I$ the 
corresponding relatively free algebra.

Let $f(y_1,\ldots, y_p,z_1,\ldots,z_q)$ be a multihomogeneous
polynomial. Then modulo the graded identity in (\ref{mainpiscalar}) we  write it as
\[
f(y_1,\ldots, y_p,z_1,\ldots,z_q) = y_1^{n_1}\ldots y_p^{n_p}
g(z_1,\ldots,z_q)
\]
where $g$ is some polynomial on the variables $z$ only. Therefore $f$
is a graded identity for $B$ if and only if $g$ is. But $g$ is a
Jordan polynomial in the variables $z$ only. Therefore $f$ is a graded
identity for $B$ if and only if $g$ is a weak Jordan identity for
the pair $(B,V)$. In this way we have to describe these weak Jordan
identities.

Define $M$ to be the subalgebra of $L=J(X)/I$ generated by the
variables $Z$.

\begin{lemma}
\label{Mgraded}
The algebra $M=M^{(0)}\oplus M^{(1)}$ is $\mathbb{Z}_2$-graded. The subalgebra
$M^{(0)}$ is spanned by all products $(z_{i_1}z_{j_1}) \ldots (z_{i_k}z_{j_k})$
while the vector space $M^{(1)}$ is spanned by all $z_{i_0}(z_{i_1}z_{j_1})
\ldots (z_{i_k}z_{j_k})$.
\end{lemma}

\pr
It is clear that the above decomposition of $M$ is a grading. The
other two statement of the lemma are equally trivial.
\epr

We denote by $C$, respectively $C_n$, the Clifford algebra of
the vector space $V$, respectively $V_n$. As we already mentioned, $C$
and $C_n$ are the associative envelopes of the special Jordan algebras
$B$ and $B_n$, respectively. The weak \textsl{associative} identities
for the pairs $(C,V)$ and $(C_n,V_n)$ were described in \cite{vdpk,
 vdpkjims} over a field of characteristic 0, and in \cite{pkfbi} over
an infinite field of characteristic different from 2. The paper
\cite{pkfbi} made use of the invariants of the orthogonal group as
described by De Concini and Procesi in \cite{dcp}. Note that the
description of these invariants in \cite{dcp} is characteristic-free.

Such a description is given in terms of double tableaux. We recall
briefly the main notions and results we shall need here. Let $t_{ij}$,
$i=1$, 2, \dots, $j=1$, 2, \dots, $n$ be commuting variables, and
consider the (formal) vectors $t_i=(t_{i1},\ldots, t_{in})$. Let
$U$ be the free $K[t_{ij}]$-module freely generated by $t_i$, $i=1$,
2, \dots{}  A nondegenerate symmetric bilinear form on $U$ is defined
as $t_i\circ t_j = t_{i1}t_{j1}+\cdots+t_{in}t_{jn}$. The description
of the polynomial algebra $R=K[(t_i\circ t_j)]$ was given in
\cite[Section 5]{dcp}.  We consider the double tableau
\begin{equation}
\label{doubletab}
T=\left(\begin{array}{cccc|cccc}
p_{11}&p_{12}&\ldots&p_{1m_1} & q_{11}&q_{12}&\ldots&q_{1m_1}\\
p_{21}&p_{22}&\ldots&p_{2m_2} &
q_{21}&q_{22}&\ldots&q_{2m_2}\\
&&&&&&&\\
\multispan 4\dotfill & \multispan 4\dotfill\\
&&&&&&&\\
p_{k1}&p_{k2}&\ldots&p_{km_k} & q_{k1}&q_{k2}&\ldots&q_{km_k}
\end{array}
\right)
\end{equation}
In it we have  $m_1\ge m_2\ge\ldots\ge m_k\ge 0$, and $p_{ij}$ and
$q_{ij}$ are integers.

Suppose $T=(p_1\, p_2\, \ldots\, p_m\mid q_1\, q_2\, \ldots\, q_m)$ is
a double row tableau filled with positive integers. We associate to it
the polynomial $\widetilde\varphi(T)\in R$:
\[
\widetilde\varphi (T)=\sum (-1)^\sigma (t_{p_1}\circ t_{q_{\sigma(1)}})(t_{p_2}\circ
t_{q_{\sigma(2)}}) \ldots (t_{p_m}\circ t_{q_{\sigma(m)}}).
\]
Here $\sigma$ runs over the symmetric group $S_m$ and
$(-1)^\sigma$ is the sign of $\sigma$. Clearly $\widetilde\varphi(T)
= \det \bigl((t_{p_i}\circ t_{q_j})\bigr)$ where $1\le i, j\le m$. If
$p_i=p_j$ or $q_i=q_j$ for some $i\ne j$ then $\widetilde\varphi(T) = 0$.

In general, if  $T^{(1)}$, $T^{(2)}$, \dots, $T^{(k)}$ are the rows of the double
tableau $T$ filled with positive integers then we associate to $T$ the polynomial
$\widetilde\varphi(T)=
\widetilde\varphi(T^{(1)})\, \widetilde\varphi(T^{(2)}) \ldots
\widetilde\varphi(T^{(k)})$.

The double tableau $T$ of the form (\ref{doubletab}) is \textsl{doubly
 standard} if the inequalities
$p_{i1}<p_{i2}<\ldots <p_{im_i}$, $q_{i1}<q_{i2}<\ldots <q_{im_i}$,
$p_{ij}\le q_{ij}$, $q_{ij}\le p_{i+1,j}$ hold for all $i$ and $j$. In
other words if we form the ordinary tableau by inserting each row of
$q_{ij}$ just below its counterpart $p_{ij}$ the resulting tableau
will be standard (that is its entries increase strictly along the
rows, and increase with possible repetitions along the columns). One
of the main results of \cite{dcp} is the following.

\begin{theorem}[\cite{dcp}, Theorem 5.1]
\label{dcp1}
The polynomials $\{\widetilde\varphi(T)\}$ where
$T$ runs over all doubly standard tableaux (\ref{doubletab}) of
positive integers,  such that $m_1\le n$,
form a basis of the vector space $R$ over $K$.
\end{theorem}

We shall call the double tableaux of the type (\ref{doubletab}) simply
\textsl{tableaux} if all their entries are positive
integers. If only $p_{11}=0$ and all remaining entries of $T$ are
positive integers we call it \textsl{0-tableau}. If $T$ is a 0-tableau
consisting of a single row we associate to it the polynomial
\[
\widetilde\varphi (T)=\sum (-1)^\sigma t_{q_{\sigma(1)}}
(t_{p_2}\circ
t_{q_{\sigma(2)}}) \ldots (t_{p_m}\circ t_{q_{\sigma(m)}}).
\]
If $T^{(1)}$, $T^{(2)}$, \dots, $T^{(k)}$ are the rows of the
0-tableau $T$ then we associate to $T$ the polynomial
$\widetilde\varphi(T)=
\widetilde\varphi(T^{(1)})\, \widetilde\varphi(T^{(2)}) \ldots
\widetilde\varphi(T^{(k)})$. Pay attention to the fact that the rows $T^{(2)}$,
\dots, $T^{(k)}$ are tableaux but not 0-tableaux.

Clearly all the above holds if we substitute $R$ by $M$ (that is $t_i$
by $z_i$). Formally
speaking one has to let $n\to\infty$ in order to justify the
statements for the infinite dimensional case but this is evidently
true.

\begin{proposition}
\label{dtspan}
The vector space $M^{(0)}$ has a basis consisting of all polynomials
associated to doubly standard tableaux. Also $M^{(1)}$ has a basis
consisting of all polynomials associated to doubly standard
0-tableaux.
\end{proposition}

\pr
The assertion for $M^{(0)}$ follows immediately from the above cited
Theorem~\ref{dcp1} of De Concini and Procesi, \cite{dcp}. The one for
$M^{(1)}$ also follows from \cite{dcp} in the following way. Our symmetric
bilinear form is nondegenerate. Let $T$ be some 0-tableau and consider
$\widetilde\varphi (T)$. Take a new variable $z_0$, then $z_0\circ
\widetilde\varphi (T)$ is represented by a double tableau, and we
apply the argument above. Then $z_0\circ \widetilde\varphi (T)$ will be
a linear combination of standard tableaux (the \textsl{straightening
 algorithm} from \cite{dcp}). But in a standard tableaux the leftmost
entry of the first row must correspond to $z_0$, and we are done.
\epr

Recall that when dealing with weak identities we consider them in the
free Jordan algebra $J(X)$.

\begin{theorem}
\label{mainthmweak}
1. The weak Jordan identities for the pair $(B,V)$ are consequences of the
polynomial $(x_1x_2,x_3,x_4)$.

2. The weak Jordan identities for the pair $(B_n,V_n)$ follow from
$(x_1x_2,x_3,x_4)$ and $f_n=\sum (-1)^\sigma x_{\sigma(1)} (x_{n+2}
x_{\sigma(2)}) \ldots (x_{2n+1}x_{\sigma(n+1)})$.
\end{theorem}

\pr
The first assertion of the theorem is a straightforward application of
Theorem~\ref{dcp1}. The same for the second assertion. (Note that the
polynomial $f_n$  ''kills'' all tableaux whose first row is of length
$\ge n+1$.)
\epr

\begin{corollary}
\label{mainresscalar}
1. The ideal of the graded identities for the Jordan algebra $B$ (with
the scalar grading) coincides with the ideal $I$ generated by the

polynomial from (\ref{mainpiscalar}).

2. The ideal of the graded identities for the Jordan algebra $B_n$
with the scalar grading is generated by (\ref{mainpiscalar}) and by
the identity
\[
g_n=\sum (-1)^\sigma z_{\sigma(1)} (z_{n+2}
z_{\sigma(2)}) \ldots (z_{2n+1}z_{\sigma(n+1)}), \quad \sigma\in S_{n+1}.
\]

3. The graded identities for the Jordan algebra of the symmetric
$2\times 2$ matrices (with the scalar grading) follow from
(\ref{mainpiscalar}) and $\sum (-1)^\sigma z_{\sigma(1)} (z_4
z_{\sigma(2)})(z_5 z_{\sigma(3)})$ where $\sigma$ runs over $S_3$.
\end{corollary}

\pr
Statement (3) is a particular case of (2). Also (1) and (2) are
immediate due to Theorem~\ref{mainthmweak} and to the remarks
preceding Lemma~\ref{Mgraded}.
\epr

\begin{center}
{\large\textbf{Acknowledgments}}
\end{center}

\noindent
We thank the Referee whose precise comments improved the
exposition of the paper. The Referee's suggestion for
including Remark~\ref{z3grading} in the final version of the
paper was taken into account with gratitude. This paper was written in
its final form while the first author was visiting the Mathematical
Institute of the University of Oxford. The first author expresses his
thanks to the colleagues of the Algebra group there for their
hospitality.

\end{document}